\documentclass{amsart}

\newtheorem{theorem}{Theorem}[section]
\newtheorem{lemma}[theorem]{Lemma}
\newtheorem{proposition}[theorem]{Proposition}

\newtheorem{example}[theorem]{Example}

\theoremstyle{definition}
\newtheorem{definition}[theorem]{Definition}

\theoremstyle{remark}

\theoremstyle{definition}

\newcommand{\equref}[1]{(\ref{#1})}

\newcommand{\pp}{{\mathbf  P}}
\newcommand{\aaf}{{\mathbf  A}}

\newcommand{\rr}{{\mathbf  R}}

\newcommand{\cc}{{\mathbf  C}}

\newcommand{\ra}{\rightarrow}

\newcommand{\st}{\,\,\mbox{such that}\,\,}

\newcommand{\lb}{{\mathcal{O}}}
\newcommand{\blowup}[1]{{\widetilde{ #1 }}}

\newcommand{\mult}{{\mathrm {mult}}}

\renewcommand{\dim}{{\mathrm {dim}}}

\newcommand{\eps}{\epsilon}
\begin{document}
\title{A lower bound on Seshadri constants of hyperplane bundles on threefolds}
\author{Kungho Chan}
%
%\institute{\textsc{Kungho Chan}\at KIAS Hoegiro 87(207-43
%Cheongnyangni 2-dong),Dongdaemun-gu,Seoul 130-722,Korea}
%
\date{Received: date / Revised version: date}
\maketitle
\begin{abstract}
We give the lower bound on Seshadri constants for the case of very
ample line bundles on threefolds. We consider the situation when the
Seshadri constant is strictly less than $2$ and give a version of
Bauer's theorem \cite[Theorem 2.1]{B1} for singular surfaces so we
can prove the same result for smooth threefolds.

% give the Mathematics Subject Classification here if any, e.g.:
% \subclass{00A00}
\end{abstract}
% If there is no abstract please say so using
% \noabstract instead
%
\section{Introduction} \label{intro}
The Seshadri constant at a given point on a smooth projective
variety was introduced by Demailly \cite{Dem} to study Fujita's
conjecture. It measures how positive a nef line bundle locally is
near a given point. Since then Seshadri constants were recognized as
interesting invariants of algebraic varieties on their own.

\begin{definition}
Suppose $X$ is a projective variety of dimension $n$ and $L$ is a
nef line bundle over $X$. Let $x$ be a point on $X$ and \[\pi :
{\blowup{X}} \ra X\] the blowup of $X$ at $x$ with the exceptional
divisor $E$. Then, we define the {\it Seshadri constant} $\eps(L,x)$
of $L$ at $x$ as
\[\eps(L,x) := \sup\{ \alpha \ge 0 ~|~ \pi^*L - \alpha E
\mbox{~is~nef~}\}.\] Or, equivalently, it can be defined by
\[\eps(L,x) := \inf\limits_{x \in C \subset X} \left\{
\frac{L.C}{\mult_x C} \right\}\] where the infimum  is taken over
all integral curves $C \subset X$ passing through $x$
\cite[5.1.5]{L1}.
\end{definition}

For the case that $X \subset \pp^N$ is a smooth integral projective
variety and $L$ is the restriction $\lb_X(1)$ on $X$ of the
hyperplane bundle of $\pp^N$, it is easy to see that
$\eps(\lb_X(1),x) \ge 1$ for every $x \in X$. Obviously, the
equality holds if there is a line in $X$ passing through $x$. For
the case of smooth surfaces, Bauer proved the following.

\begin{theorem} \cite[Theorem 2.1]{B1} \label{bauer:thm}
Let $X \subset \pp^N$ be a smooth irreducible surface.
\begin{enumerate}
\item[(a)] $\eps(\lb_X(1),x) = 1$ if and only if $X$ contains a line
passing through $x$.
\item[(b)] Suppose $X$ is of degree $d~(\ge 4)$ and $x$ is a point on $X$. If $X$ contains no
line passing through $x$, then

\[ \eps(\lb_X(1),x) \ge \frac{d}{d-1}.\]

\item[(c)] If $X$ is of degree $d~(\ge 4)$ and $x \in X$ is a point such that the Seshadri constant
$\eps(\lb_X(1),x)$ satisfies the inequalities $1 < \eps(\lb_X(1),x)
< 2$, then it is of the form

\[\eps(\lb_X(1),x) = \frac{a}{b},\] where $a,b$ are integers with $3
\le a \le d$ and $a/2 < b < a$.

\item[(d)] All rational numbers $a/b$ with $3
\le a \le d$ and $a/2 < b < a$ occur as local Seshadri constants of
smooth irreducible surfaces in $\pp^3$ of degree $d$.
\end{enumerate}
\end{theorem}

Bauer's approach was to consider the intersection of $X$ and the
tangent plane $T_x X$ of $X$ at $x$. If the Seshadri constant
$\eps(\lb_X(1),x) < 2$, then it must be computed by a component of
the intersection and the multiplicity of the component at $x$ is
bounded above by the degree $d$ of $X$. In this paper, we prove a
similar result for any integral surfaces in $\pp^3$.

\begin{proposition} \label{surface:case}  Let $X$ be a projective
integral surface of degree $d~( \ge 3)$ in $\pp^3$. If $x \in X$ is
a point of multiplicity $m$ and $X$ contains no line passing through
$x$, then
\item[(a)]
\[\eps(\lb_X(1),x) \ge \frac{d}{d-1}.\]
\item[(b)]
If $\eps(\lb_X(1),x) < \frac{m+1}{m}$, then

\[\eps(\lb_X(1),x) = \frac{a}{b}\] for some integers $a,b$ such
that $3 \le a \le md$ and $\frac{ma}{m+1} < b \le \frac{a(d-1)}{d}
$.
\end{proposition}

By constructing a singular surface in $\pp^3$, we also show that the
lower bound $\frac{d}{d-1}$ here is optimal. Furthermore, we give
the same lower bound for smooth threefolds in projective spaces.

\begin{theorem} \label{very:ample:3fold}
Let $X$ be an irreducible smooth projective threefold of degree
$d~(\ge 4)$ in $\pp^N$, $x$ a point in $X$ and $T_x X$ the tangent
linear subspace of $\pp^N$ to $X$ at $x$. If $X$ contains no line
through $x$, then
\begin{enumerate}
\item[(a)] If $\dim(T_x X \cap X) = 0$, then \[\eps(\lb_X(1),x) \ge 2.\]
\item[(b)] If $\dim(T_x X \cap X) = 1$, then \[\eps(\lb_X(1),x) \ge \frac{d}{d-1}.\]
If furthermore $\eps(\lb_X(1),x) < 2$, then $\eps(\lb_X(1),x) =
\frac{a}{b}$ for some integers $a,b$ such that $3 \le a \le d$ and
$\frac{a}{2} < b < a$.
\item[(c)] If $\dim(T_x X \cap X) = 2$, then \[\eps(\lb_X(1),x) \ge \frac{d}{d-1}. \]
If furthermore $\eps(\lb_X(1),x) < \frac{d-1}{d-2}$, then
$\eps(\lb_X(1),x) = \frac{a}{b}$ for some integers $a,b$ such that
$3 \le a \le (d-2)d$ and $\frac{(d-2)a}{d-1} < b \le
\frac{a(d-1)}{d}$.
\end{enumerate}
\end{theorem}

Finally, we construct a threefold with finitely many singular points
in $\pp^4$ whose Seshadri constant of hyperplane bundle at a smooth
point is $\frac{d}{d-1}$.

{\it{Acknowledgement.}} {This paper is part of my PhD thesis. I
thank my advisor L. Ein for getting me interested in this question
and helping me with comments and suggestions.}

\section{Seshadri constants of hyperplane line bundles on
surfaces}

In order to prove \ref{surface:case}, we need some preparations.

\begin{definition}
Let X be a reduced irreducible projective variety of dimension $n$
in $\pp^N$ and $x$ a point on $X$. Denote $\lb_X(1) =
\lb_{\pp^N}(1)|_{X}$ by $H$, the blowup of $\pp^N$ at $x$ by
${\widetilde{\pp^N}}$, the proper transform of $X$ in
${\widetilde{\pp^N}}$ by ${\widetilde{X}}$ and the exceptional
divisor in ${\widetilde{X}}$ by $E$. The projection $\pp^N
\backslash \{x\} \ra \pp^{N-1}$ induces a morphism
${\widetilde{\pp^N}} \ra \pp^{N-1}$ such that the proper transform
of every line in $\pp^{N}$ passing through $x$ is mapped to a point
in $\pp^{N-1}$. We call the restriction on ${\widetilde{X}}$ of this
morphism the {\it inner projection} of $X \subset \pp^N$ from $x$,
that is

\[p : {\widetilde{X}} \ra \pp^{N-1}.\] If $\pi:{\widetilde{X}} \ra X$ is the blowup morphism of $X$ at
$x$, then from the construction $p^*\lb(1) = \pi ^*H - E$.
\end{definition}
It is obvious that if there is no line in $X$
passing through $x$, then $p$ is finite.

\begin{lemma} \label{no:line:lemma} With the above settings, there is no
line in $X$ passing through $x$ if and only if $\eps(H,x) > 1$.
\end{lemma}
\begin{proof}  The "$\Leftarrow$" is obvious. For "$\Rightarrow$", we assume
the hypothesis, then $p$ is finite. This implies that $p^*\lb(1)$ is
ample and hence $\pi^*H - E$ is ample. So $\eps(H,x) > 1$.
\end{proof}

\begin{lemma} \label{multiplicity:lemma} With the above settings, if  $\eps(H,x) > 1$, then
\[\deg X - r \ge \mult_x X\] where $r =
\deg p(\blowup{X})$. In particular, $\deg X - 1 \ge \mult_x X$.
\end{lemma}
\begin{proof}  By \ref{no:line:lemma}, $X$ contains no line passing
through $x$, $p$ is finite and hence
\begin{equation} \label{mult:x}
\begin{split}
\deg X - \mult_x X  &= \deg X + (-1)^n E^n \\
                    &= (\pi^*H - E)^n\\
                    &= p^*\lb(1)^n\\
                    &= \deg(p)\cdot
(\lb(1)|_{p({\widetilde{X}})})^n\\
                    &= \deg(p) \cdot \deg p(\blowup{X})\\
                    & \ge r.\\
\end{split}
\end{equation}
\end{proof}

The following two lemmas gives us a way to bound the multiplicity of
a variety at a point by its degree.

\begin{lemma} \label{on:cone}
Suppose $T \subset \pp^N$ is an integral cone of dimension $N-1$
over an integral hypersurface $B$ in $\pp^{N-1}$ with vertex at $p
\in \pp^N$. If $X$ is an integral subvariety in $T$ of codimension
one containing $p$ and $X$ contains no line passing through $p$,
then
\[ \deg X - \deg B \ge \mult_p X.\]
\end{lemma}

\begin{proof}
Since the image $p({\widetilde{X}}) = B$ and $\eps(\lb_X(1), p) >
1$, then from \equref{mult:x}, we have \[\deg X - \mult_p X =
p^*\lb(1)^{N-2}.{\widetilde{X}} = \deg(p) \cdot (\lb(1)|_B)^{N-2} =
\deg{p} \cdot \deg B \ge \deg B.\]
\end{proof}

\begin{lemma} \label{hypersurface:mult}
Suppose $X$ is a hypersurface of degree $d$ in $\pp^{n+1}$ and $x$
is a point of multiplicity $m$ on $X$. If $X$ contains no line
passing through $x$, then
\[ m \le d -n.\]
\end{lemma}

\begin{proof}
Suppose $X$ is defined by the homogeneous polynomial $f$ of degree
$d$. By choosing an appropriate coordinate system, we can assume $x
= [0:0:\cdots:0:1]$. Then

\[f = x_{n+1}^{d-1}\cdot f_1 + x_{n+1}^{d-2}\cdot f_2 + \cdots +
x_{n+1}^{d-i} \cdot f_i + \cdots + f_d\]

where $f_i$ is a homogeneous polynomial of degree $i$ in
$x_0,x_1,\dots,x_n$. Since $\mult_x X = m$, then $f_m \not = 0$ and
$f_i = 0$ for $i < m$. Assume $m > d -n$, then there are at most $n$
nonzero $f_i$'s with $d \ge i \ge m$. Then the intersection of the
cones in $\pp^{n+1}$ defined by nonzero $f_i$'s is nonempty and
contains line. This implies $X$ contains a line passing through $x$.
This is a contradiction.
\end{proof}

Now we come up to prove \ref{surface:case}.

\begin{proof}[Proof of \ref{surface:case}]  Let $TC_x(X) \subset \pp^3$ be the tangent cone of $X$ at
$x$. Suppose the local equation of $X$ around $x$ is $f$. By
choosing an appropriate coordinate system, we can assume $x =
(0,0,0)$ and $f = f_m + f_{\ge m+1}$ where $f_i$ is a homogeneous
polynomial of degree $i$. Then the local equation of $TC_x(X)$ is
$f_m$, and hence $\deg TC_x(X) = m$ and $\mult_x (TC_x(X)) = m$.
Since $TC_x(X)$ is spanned by all the tangent lines to $X$ at $x$
and $X$ contains no line through $x$, $TC_x(X)$ must intersect $X$
properly. Denote $TC_x(X) \cap X$ by $D$. Then $D \in |\lb_X(m)|$.

If $\eps(\lb_X(1),x) \ge \frac{m+1}{m}$, then $\eps(\lb_X(1),x) \ge
\frac{d}{d-1}$ since $\frac{x+1}{x}$ is decreasing for $x > 0$ and
$d-1 \ge m$ by \ref{multiplicity:lemma}. Hence it reduces to assume
\begin{equation} \label{the:critical:upperbound}
\eps(\lb_X(1),x) < \frac{m+1}{m}.
\end{equation}
This implies that there exists an integral curve $x \in C \subset X$
such that $\eps(\lb_X(1),x) \le \frac{\deg C}{\mult_xC} <
\frac{m+1}{m}$. We claim that $C$ must be an irreducible component
of $D$.  To this end, let us consider the blowup of $\pp^3$ at $x$
with the exceptional divisor $E$,
\[ \pi : {\widetilde{\pp^3}} \ra \pp^3.\] Let us set $T = TC_x(X)$.
Then,
\[\pi^*T = {\widetilde{T}} + mE\] and
\[\pi^*X = {\widetilde{X}} + mE\] where ${\blowup{T}}$ and
$\blowup{X}$ are respectively the proper transforms of $T$ and $X$.
Since $\blowup{X}|_E = \blowup{T}|_E$, then

\[(\pi|_{\blowup{X}})^*D = \blowup{D} + \mbox{the common components of }(\blowup{X}|_E)\mbox{ and }(\blowup{T}|_E) +
m \cdot \blowup{X}|_E = \blowup{D} + (m+1)(\blowup{X}|_E)\] where
$\blowup{D}$ is the proper transform of $D$.

Let $\blowup{C}$ be the proper transform of $C$.

Assume $C$ is not an irreducible component of $D$, then
\begin{equation*}
\begin{split}
m\deg C = D.C &= ((\pi|_{\blowup{X}})^*D).\blowup{C}\\
                &= (\blowup{D} + (m+1)(\blowup{X}|_E)).\blowup{C}\\
                &\ge (m+1) \mult_x C.\\
\end{split}
\end{equation*}
This is a contradiction and thus $C$ is an irreducible component of
$D$.

Suppose the tangent cone $T = \sum a_i T_i$ where $T_i$'s are
integral cones in $\pp^3$ with vertices at $x$ and $\sum a_i \deg
T_i = m$. Then, $C \subset (T_i \cap X)$ for some $i$ and hence
\[\deg C \le \deg X
\cdot \deg T_i = d \deg T_i.\] Moreover,  by \ref{on:cone}
\begin{equation*} \label{seshadri:curve:mult}
\frac{\deg C}{\mult_x C} \ge \frac{\deg C}{\deg C - \deg T_i} \ge
\frac{d \deg T_i}{ d \deg T_i - \deg T_i} =  \frac{d}{d-1}.
\end{equation*}
So we have $\eps(\lb_X(1),x) \ge \frac{d}{d-1}$. This is for part
(a).

For part (b), if \equref{the:critical:upperbound} is satisfied,
$\eps(\lb_X(1),x)$ must be computed by one of the irreducible
components of $D$. Let $C$ be such a component and contained in an
integral cone $T_i \subset T$ then \[\eps(\lb_X(1),x) = \frac{\deg
C}{\mult_x C}.\] Since $X$ contains no line passing through $x$, $1
< \deg C (\le \deg D = md)$ and $(\frac{m\deg C}{m+1} <) \mult_x C
\le \frac{(d-1) \deg C}{d}$. Since $\eps(\lb_X(1),x) < \frac{m+1}{m}
\le 2$, the $\deg C$ cannot equal to $2$ and hence $\deg C \ge 3$.
\end{proof}

The following example indicates that the lower bound $\frac{d}{d-1}$
is optimal.

\begin{example} \label{surface:example}
For any $d \ge 4$ and $2 \le m \le d-2$, there exists an integral
surface $S \subset \pp^3$ containing no line through a point $0 \in
S$ such that

\[ \mult_0 S = m \mbox{~and~} \eps( \lb (1),0) =
\frac{d}{d-1}.\]
\end{example}

\begin{proof}[Construction of \ref{surface:example}]
Roughly speaking, we construct $S$ by prescribing its tangent cone
$T$, and then we have to show the irreducibility of $S$ and it
contains no line through a particular point. Moreover, we also have
to show the irreducibility of the intersection $S \cap T$.

Choose a general integral plane curve $B$ of degree $m$ and let $T$
be the cone over $B$ with vertex $0 \in \pp^3$. Let the blowup of
$\pp^3$ at $0$ be

\[\pi : \blowup{\pp^3} \ra \pp^3\] with the exceptional divisor $E$ and $\blowup{T}$ the proper
transform of $T$. Consider on $\blowup{\pp^3}$ the line bundle
\[ L := \pi^*\lb_{\pp^3}(d) - (d-1)E.\] It is clear that $L$ is
globally generated.

Let $V$ be the image of the global section restriction morphism,

\[V := {\mathrm{Image}}\left(H^0(\blowup{\pp^3},L) \ra
H^0(\blowup{T},L|_{\blowup{T}})\right).\]

Take a general effective divisor $D \in |L|$. Then $D|_E \not\supset
B = \blowup{T}|_E$. Set $Y := \pi_*D$ and $C := Y|_T$, then \[ \deg
C = dm.\] Moreover, set $\sigma := \pi|_{\blowup{T}} : \blowup{T}
\ra T$, then

\[\sigma^*(C) = (D + (d-1)E)|_{\blowup{T}} = \blowup{C} +
(d-1)B\] where the $\blowup{C}$ is the proper transform of $C$.

Thus, \[\mult_0 (C) = \blowup{C}.B = (\sigma^*(C) - (d-1)B).B =
(d-1)m.\]

Now we show the irreducibility of $C$. Consider the morphism induced
by $|V|$,

\[ \phi : \blowup{T} \ra \pp^{\dim |V|}.\]

Since $(L|_{\blowup{T}})^2 = (d^2- (d-1)^2)m > 0$, the image of
$\phi$ cannot be $1$-dimensional. Then, by the Bertini theorem for
general linear sections \cite[3.3.1]{L1}, $\blowup{C}$ is
irreducible and so is $C$.

With the $T$, $Y$ and $C$ constructed above, we can construct the
singular surface $S$.

Let $U \ni 0$ be the complement of a hyperplane in $\pp^3$, which is
affine $\aaf^3$.

By choosing an appropriate coordinate system, we can assume the
defining equation of $T$ in $\aaf^3$ is $f_{m} =0$ with the vertex
$(0,0,0)$ and the defining equation of $Y$ is $f_{d-1} + f_d=0$
where $f_i$ is a homogeneous polynomial of degree $i$ in the
variables $x_1,x_2$ and $x_3$.

Now, we define the surface $S|_U$ by the equation

\[f = f_{m} + f_{d-1} + f_{d} = 0.\] Clearly,

\[\mbox{~the tangent cone~}TC_0(S) = T \mbox{~and~} S \cap T = C.\]

Since $T$ is irreducible, then only one irreducible component of $S$
can pass through the point $0$. If $S$ is reducible, then $\deg C <
dm$. This is impossible when we choose $C$ with $\deg C = dm$. Thus,
$S$ is irreducible and

\[ \frac{\deg C}{\mult_0 C} = \frac{dm}{(d-1)m} =
\frac{d}{d-1} < \frac{m+1}{m}.\] From the proof of
\ref{surface:case}, $C$ computes the Seshadri constant

\[ \eps(\lb(1),0) = \frac{d}{d-1}.\]

To complete the proof, it remains to show that $S$ contains no line
through $0$. Assume to the contrary that $S$ contains a line through
$0$. Let the line be $\ell := \{(at,bt,ct)| t \in \cc\}$. Then,
\[f(at,bt,ct) = 0 ~~\forall~ t.\] This implies that $\ell \subset (T \cap C(f_{d-1}) \cap C(f_{d}))$
where $C(f_{j})$ is the cone in $\aaf^3$ defined by the homogeneous
polynomial $f_j$. However, this situation can be eliminated when we
made a general choice of $D$ above.
\end{proof}

\section{Seshadri constants of hyperplane bundles on smooth
threefolds}

To extend Bauer's approach to smooth threefold $X$, one has to
consider the intersection of the $3$-dimensional tangent space $T_x
X$ and $X$. The possible components of the intersection are points,
curves and surfaces. It turns out that \ref{surface:case} can deal
with the surface components.

\begin{proof}[Proof of \ref{very:ample:3fold}]
We start from a simple argument. Suppose $C \subset X$ is an
integral curve through $x$. If
\begin{equation} \label{key:criterion}
C \not\subset T_x X \cap X,
\end{equation} one can take  a
hyperplane $H \subset \pp^N \st H \supset T_x X$ and $H \not\supset
C$. Then \[\deg C = (H|_X).C \ge 2 \cdot \mult_x C.\]

For (a), \equref{key:criterion} is always true for every integral
curve through $x$. So, $\eps(\lb_X(1),x) \ge 2$.

For (b), assume $\frac{\deg C}{\mult_x C} < 2$. Then $C$ must be an
irreducible component of $T_x X \cap X$. Otherwise,
\equref{key:criterion} is true for $C$. Now $C \subset T_x X \cap
X$, then $\deg C \le d$. By \ref{multiplicity:lemma}
\[ \frac{\deg C}{\mult_x C} \ge \frac{\deg C}{\deg C-1} \ge
\frac{d}{d-1}.\] This induces $\eps(\lb_X(1),x) \ge \frac{d}{d-1}$.

In this situation, $\eps(\lb_X(1),x)$ must be computed by one of the
integral curves in $T_x X \cap X$ with $2 < \deg C \le d$ and
$\frac{\deg C}{2} < \mult_x C < \deg C$.

For (c), assume $\frac{\deg C}{\mult_x C} < 2$. We can assume that
there is an irreducible and reduced surface $S \subset T_x X \cap X
\st C \subset S$. If not,  either \equref{key:criterion} is true for
$C$ or $C$ is an irreducible component of $T_x X \cap X$. For the
latter situation $\deg C \le d$, we will have $\frac{\deg C}{\mult_x
C} \ge \frac{d}{d-1}$ as the situation in (b).

Now $S$ is a projective surface in $T_x X (\cong \pp^3)$ and $\deg S
\le d$. By \ref{surface:case},

\begin{equation} \label{seshadri:curve:mult2}
\frac{\deg C}{\mult_x C} \ge \eps(\lb_X(1)|_{S},x) \ge
\frac{d}{d-1}.
\end{equation} Combining both cases, we have $\eps(\lb_X(1),x) \ge \frac{d}{d-1}$.

Assume furthermore $\eps(\lb_X(1),x) < \frac{d-1}{d-2}$.

Suppose $\mult_x S = m$. Since $d-2 \ge m$ from
\ref{hypersurface:mult}, \[\eps(\lb_X(1),x) < \frac{d-1}{d-2} \le
\frac{m+1}{m}.\] Hence, as the part (b) of \ref{surface:case}, we
have
\[\eps(\lb_X(1),x) = \frac{a}{b}\] for some integers $a,b$ such
that $3 \le a \le (d-2)d$ and $\frac{(d-2)a}{d-1} < b \le
\frac{a(d-1)}{d}$.
\end{proof}

This result implies that if $\eps(\lb_X(1),x) < \frac{d-1}{d-2}$,
then the number of possible values of $\eps(\lb_X(1),x)$ is finite.
Indeed, this phenomenon was generally proved for surfaces. In
\cite{Og}, Oguiso proved that if $(f:{\mathcal X} \ra {\mathcal B},
{\mathcal L})$ is a family of polarized surfaces of degree $d$ and
$\Sigma := \{\eps({\mathcal L}_t,x_t)~|~t\in {\mathcal B},~x_t \in
{\mathcal X}_t\}$, then for each given number $\alpha \in \rr$ such
that $\alpha < \sqrt{d}$, the set $\Sigma \cap (0, \alpha]$ is
finite.

Here we give an example of a $3$-dimensional integral hypersurface
$X$ of degree $d \ge 4$ in $\pp^4$, which is smooth at a point $0$
and it computes $\eps(\lb_X(1),0) = \frac{d}{d-1}$.

\begin{example} \label{3fold:example}
For any $d \ge 4$, there exists an integral hypersurface $X \subset
\pp^4$ containing no line through a point $0 \in X$ and smooth at
$0$ such that

\[\eps( \lb_X(1),0) =
\frac{d}{d-1}.\]
\end{example}

\begin{proof}[Construction of \ref{3fold:example}] The construction
is a straight forward extension of Bauer's example \cite[Lemma
2.2]{B1}.

We start with a surface $S$ from the example \ref{surface:example}.
Let the defining equation of $S|_U$ is:

\[g(x_1,x_2,x_3) = 0\] where $U \cong \aaf^3$ is the affine open neighborhood containing the point $0$. Choose a general smooth integral surface $S'$ in $\pp^3$ of degree $d-1$ satisfying

\begin{enumerate}
\item[(1)] $S'$ does not contain the point $0=(0,0,0)$.
\item[(2)] $S'$ does not contain any $1$-dimensional singular locus component of $S$.
\item[(3)] If a point $p$ is in the intersection of $S$ and $S'$, and if $S$ is smooth at $p$, then the tangent planes $T_p(S')$ and $T_p(S)$ are not equal.
\end{enumerate}

Let the defining equation of $S'$ be
\[g'(x_1,x_2,x_3) = 0.\] Then, we define a hypersurface $X$ in the affine open neighborhood $\aaf^4 \supset U$ by the equation
\[f = g + x_4g' = 0.\]

We claim that $X$ is smooth at the point $0=(0,0,0,0)$ and $X$ is
smooth except perhaps finitely many points.

Suppose $X$ is singular at the point $p = (p_1,p_2,p_3,p_4)$.  The
partial derivative of $f$ with respect to $x_4$ gives
$g'(p_1,p_2,p_3) = 0$, and hence $f(p) = 0$ implies $g(p_1,p_2,p_3)
=0$. Therefore, the point $p' = (p_1,p_2,p_3)$ is an intersection
point of $S$ and $S'$. By the condition (1), $X$ must be smooth at
the point $0$.

Let $j^1_p(h)$ be the linear term of $\theta(h)$ where $\theta$ is a
linear translation of the point $0$ to the point $p =
(p_1,p_2,p_3,p_4)$ where $h$ is a polynomial in $x_i$ for $1\le i
\le 4$.

Since $X$ is singular at $p$, then $j^1_p(f) = 0$. On the other
side,

\[0= j^1_p(f) = j^1_{p'}(g) + p_4 \cdot j^1_{p'}(g').\]

If $p_4 \not= 0$ and $j^1_{p'}(g) =0$, then \[j^1_{p'}(g') = 0.\]
This is impossible when $S'$ is smooth.

If $p_4 \not= 0$ and $j^1_{p'}(g) \not= 0$, then \[T_{p'}(S') =
T_{p'}(S).\] This cannot be true by the condition (3).

If $p_4 = 0$, then $j^1_{p'}(g) = 0$. It means $S$ is singular at
$p'$. But only finitely many intersection points of $S$ and $S'$ can
be in this situation by the condition (2). This completes the proof
of the claim.

Since $S$ is singular at $0$, it is clear that the hyperplane $\{x_4
=0\} = T_0 X$ and hence $S = X \cap T_0 X$.

Since $X$ is singular at most finitely many points, then $X$ is
irreducible.

Moreover, $X$ does not contain a line through $0$. To this end,
assume to the contrary that $X$ contains a line through $0$. Let the
line be $\ell := \{(a_1t,a_2t,a_3t,a_4t)| t \in \cc\}$. Then,
\[f(a_1t,a_2t,a_3t,a_4t) = 0 ~~\forall~ t.\] However , if $a_4 \not= 0$,
then it is impossible when $g'$ has nonzero constant term, and $g$
has no constant term and no linear term. If $a_4 = 0$, then the line
$\ell \subset S$. This is not possible when $S$ contains no line
through $0$. Therefore, we have
\[\eps(\lb_X(1),0) = \eps(\lb_S(1),0) = \frac{d}{d-1}.\]
\end{proof}

%
% BibTeX users please use
\bibliographystyle{alpha}
\bibliography{references}

\begin{thebibliography}{Dem92}

\bibitem[Bau99]{B1}
Thomas Bauer.
\newblock Seshadri constants on algebraic surfaces.
\newblock {\em Math. Ann.}, 313(3):547--583, 1999.

\bibitem[Dem92]{Dem}
Jean-Pierre Demailly.
\newblock Singular {H}ermitian metrics on positive line bundles.
\newblock In {\em Complex algebraic varieties (Bayreuth, 1990)}, volume 1507 of
  {\em Lecture Notes in Math.}, pages 87--104. Springer, Berlin, 1992.

\bibitem[Laz04]{L1}
Robert Lazarsfeld.
\newblock {\em Positivity in algebraic geometry. {I}}, volume~48 of {\em
  [Results in Mathematics and Related Areas. 3rd Series. A Series of Modern
  Surveys in Mathematics]}.
\newblock Springer-Verlag, Berlin, 2004.
\newblock Classical setting: line bundles and linear series.

\bibitem[Ogu02]{Og}
Keiji Oguiso.
\newblock Seshadri constants in a family of surfaces.
\newblock {\em Math. Ann.}, 323(4):625--631, 2002.

\end{thebibliography}
\end{document}